\documentclass[10pt]{article}
\usepackage[left=2.54cm,top=2.54cm,right=2.54cm,bottom=2.54cm]{geometry}
\usepackage{fancyhdr}
\usepackage{afterpage}
\usepackage{setspace}
\newfont{\toto}{msbm10 at 12 pt}
\newfont{\ithd}{cmr9}

\usepackage{booktabs}
\usepackage{amsmath,amsthm,amsfonts,amssymb}
\usepackage[pdftex]{graphicx}
\usepackage{epstopdf}
\usepackage[T1]{fontenc}
\usepackage{float}
\usepackage{placeins}
\usepackage{url}
\usepackage{hyperref}
\title{
\bf Chaotic Behavior of Stiff ODEs and Their Derivatives: An Illustrative Example 
}
\author{
 Emre \"{O}zkaya$^{1}$, Nicolas R. Gauger$^{1}$, Anil Nemili$^{2}$\\\\
$^{1}$ Chair for Scientific Computing, TU Kaiserslautern, 67663 Kaiserslautern, Germany\\
$^{2}$ Department of Mathematics, BITS Pilani, Hyderabad Campus, Hyderabad 500078, India\\
}
\date{}

\begin{document}

\maketitle
\afterpage{\fancyhead{}}


In the following, an illustrative example concerning difficulties in differentiating stiff ODEs is presented. In the given example, the solution of a completely deterministic system becomes chaotic due to computational noise introduced by the numerical algorithm.

Use of certain numerical algorithms with stiff ODEs can lead to numerically unstable or noisy solutions. Such an example is given by Moler in a Mathworks article \cite{Moler}: 
\begin{equation*}
\dot x=\alpha x^2-\beta x^3,
\end{equation*} 
with $x(0)=\sigma$ being the initial condition. This ODE is a model of flame propagation. The variable $x$ represents the radius of a ball shaped flame front, which grows rapidly in a process, for example when a match is lighted. Starting from a very small surface ($\sigma$), the flame front grows until it reaches a critical size. Then it remains at that size because the amount of oxygen being consumed by the combustion in the interior of the ball balances the amount of oxygen available through the surface. The parameters $\alpha$ and $\beta$ are scalar quantities, which model the effect of surface area and volume occupied by the flame front. Note that this ODE has steady-state solutions for certain values of $\alpha$ and $\beta$.

Combined with a certain integration method, the solution of the ODE has a very stiff behavior if $\sigma$ is small. In the following, we have taken $\sigma=0.0001$ for all numerical tests. As far as the time integration method is concerned, we have taken the explicit Euler method. Trying different $\Delta t$ values, it is possible to observe four different behavior of the ODE solution:
\begin{itemize}
\item CASE I: Monotone solution \\
\item CASE II: Oscillating solution \\
\item CASE III: Chaotic solution \\
\item CASE IV: Divergent solution 
\end{itemize}

For each case, the derivative results are computed using the tangent linear solver that is generated by the forward mode of AD. Since the forward mode of AD corresponds to exact linearization of an arbitrary function, it is a very suitable method to investigate the effect of differentiation for different values of $\Delta t$. For the numerical experiments, we consider the scalar objective function
\begin{equation*}
J=x_{ave}=\frac{1}{100000-2000}\sum_{t_n=2000}^{100000} x(t_n), 
\end{equation*}

which is the time average of the solution $x$ over a long time history. To be able get rid of transient effects, the averaging process starts after $2000$ time iterations. For a specified objective function, the tangent linear solver computes the derivative of $J$ with respect to the ODE model parameters $\alpha$ and $\beta$. For the sake of simplicity, we focus on the derivatives $d J/ d \alpha$ in the following.   

As one might expect, the most trivial case is the divergent solution. If the $\Delta t$ is taken too large, the solution explodes quickly after a few iterations. Similarly, the derivative values also explode. This is no surprise since explicit Euler scheme is well known to be unstable above a certain $\Delta t$.

\begin{figure}
  \centering
    \includegraphics[trim=1.5cm 0cm 0.5cm 0cm, clip=true, scale=0.7]{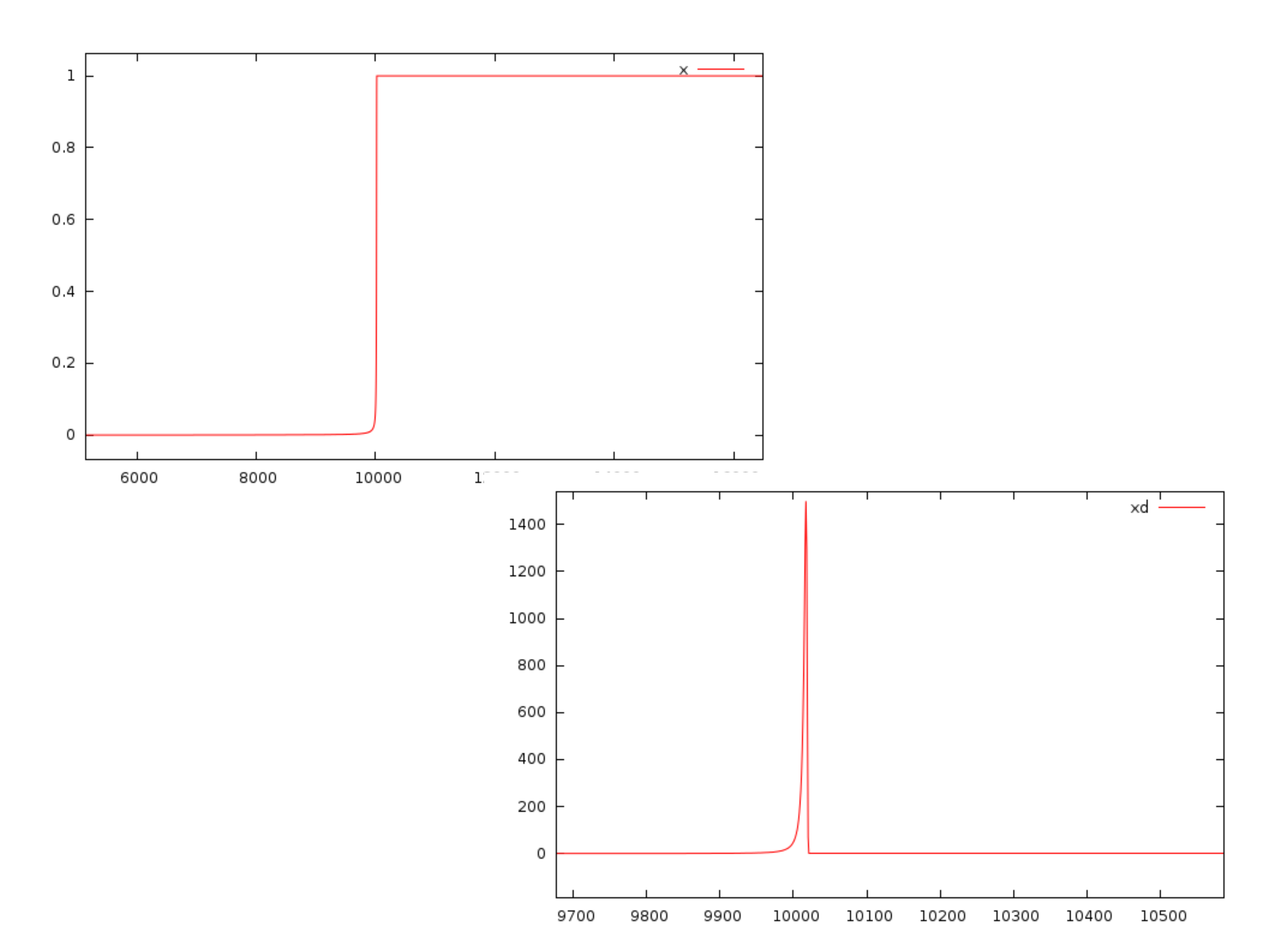}
    \caption{$x$ and $xd=dx/d\alpha$ for $\Delta t=1.0$, $\alpha=1.0$ and $\beta=1.0$ }
    \label{fig1}
\end{figure}

The case I is also simple. If $\Delta t$ is chosen small enough, the solution of the ODE is monotone and smooth. In the Fig. \ref{fig1}, the ODE solution as well as the derivative of $x$ w.r.t. $\alpha$, which is denoted by $xd$, are shown. $\Delta t$ for this case is taken as $1.0$.

\begin{figure}
  \centering
    \includegraphics[scale=0.7]{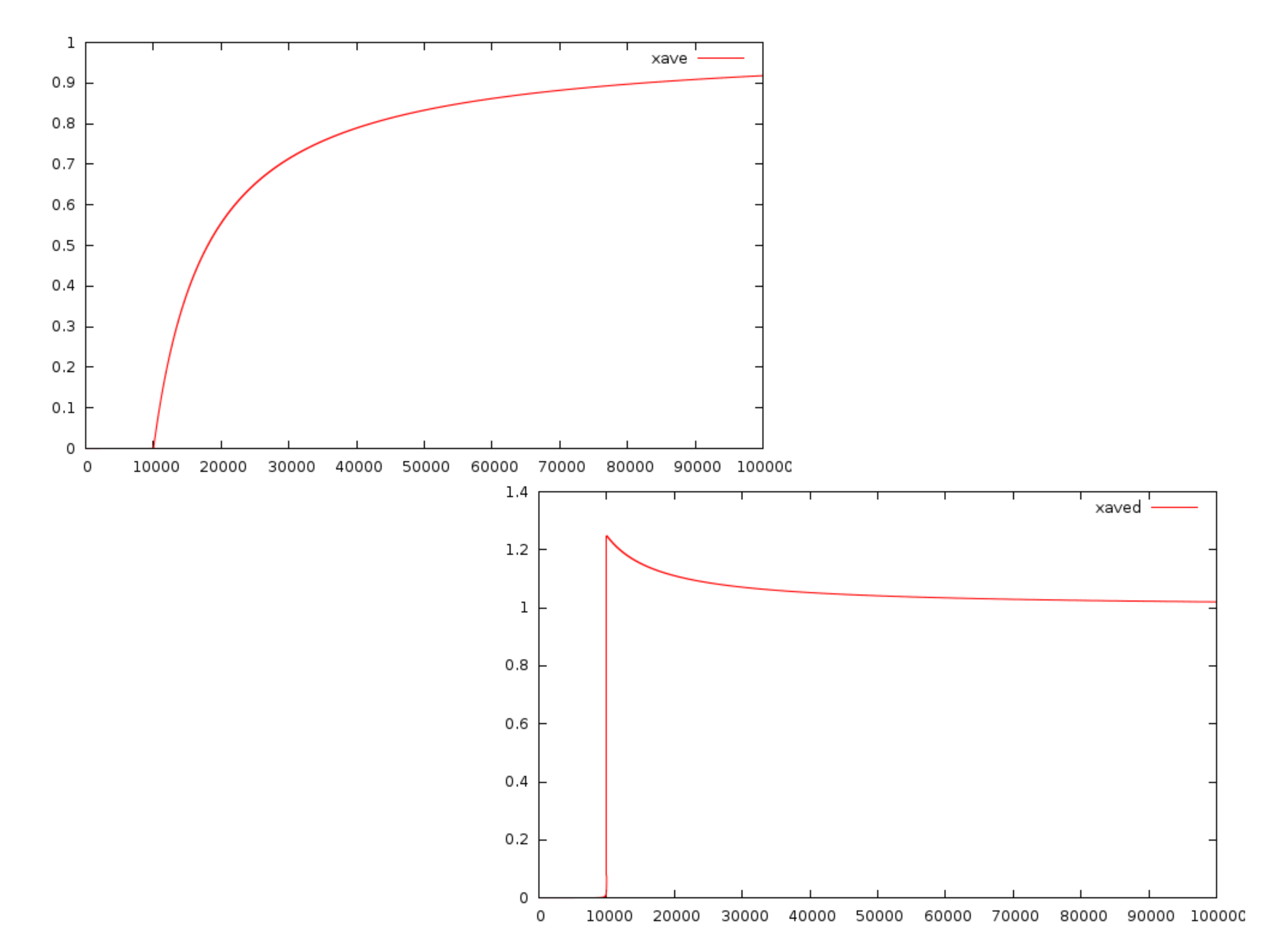}
    \caption{$J$ and $dJ/d\alpha$ for $\Delta t=1.0$, $\alpha=1.0$ and $\beta=1.0$ }
    \label{fig2}
\end{figure}

Initially, the variable $x$ grows very slowly. Around $t=10000$, however, it rapidly grows and settles to the steady-state value of $x=1.0$. The derivative value $xd$ also grows very slowly, and then jumps suddenly to a maximum value of $1499.79$ around $t=10017$. Again very rapidly it drops to $1.0$ and settles down. If we look at the convergence histories of the objective function $J=x_{ave}$ and its derivative w.r.t. $\alpha$, we again see a smooth behavior. These results can be seen in Fig. \ref{fig2}.

\begin{figure}
  \centering
    \includegraphics[scale=0.7]{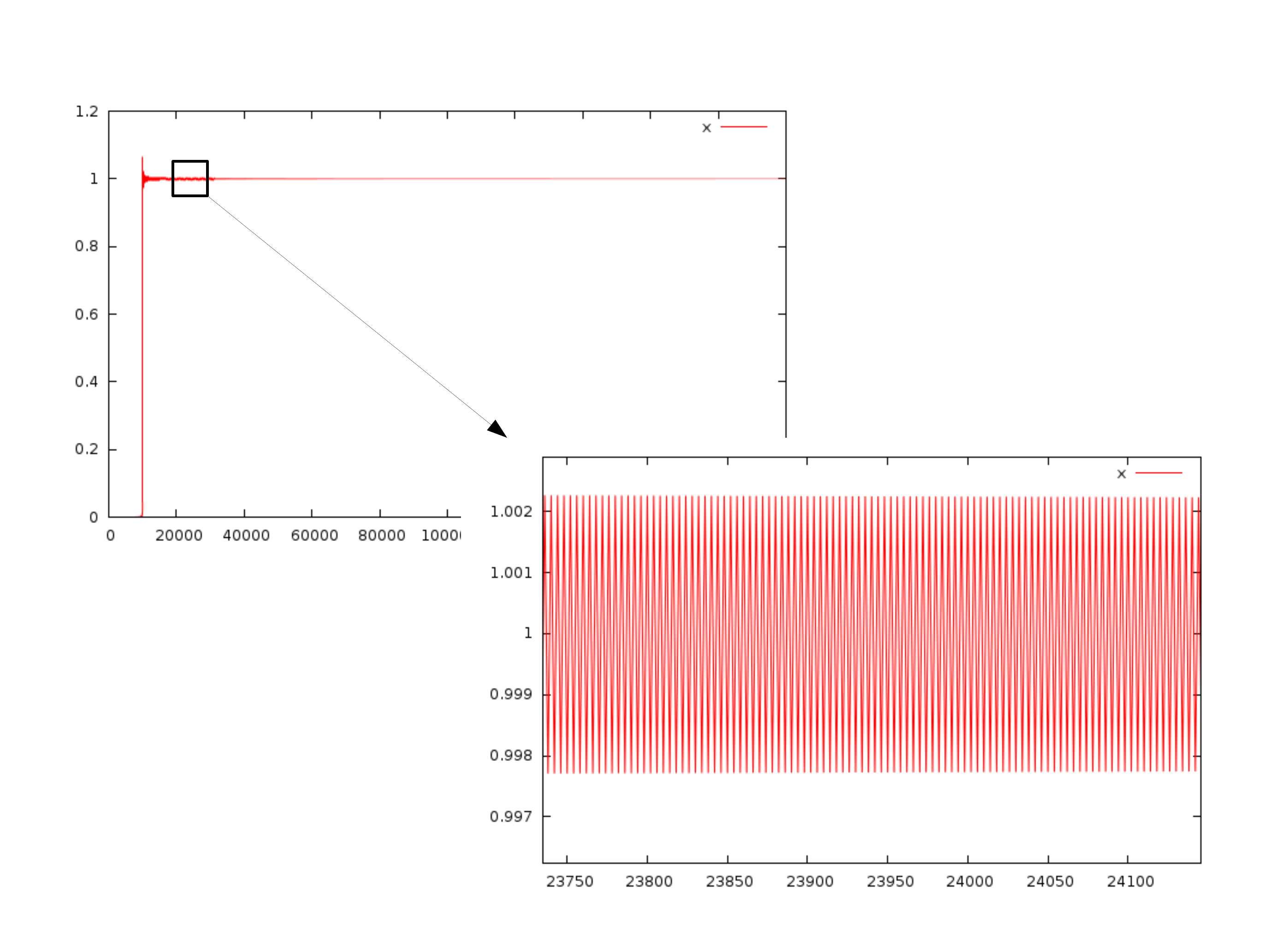}
    \caption{$x$ for $\Delta t=2.0$, $\alpha=1.0$ and $\beta=1.0$ }
    \label{fig3}
\end{figure}

\begin{figure}
  \centering
    \includegraphics[scale=0.7]{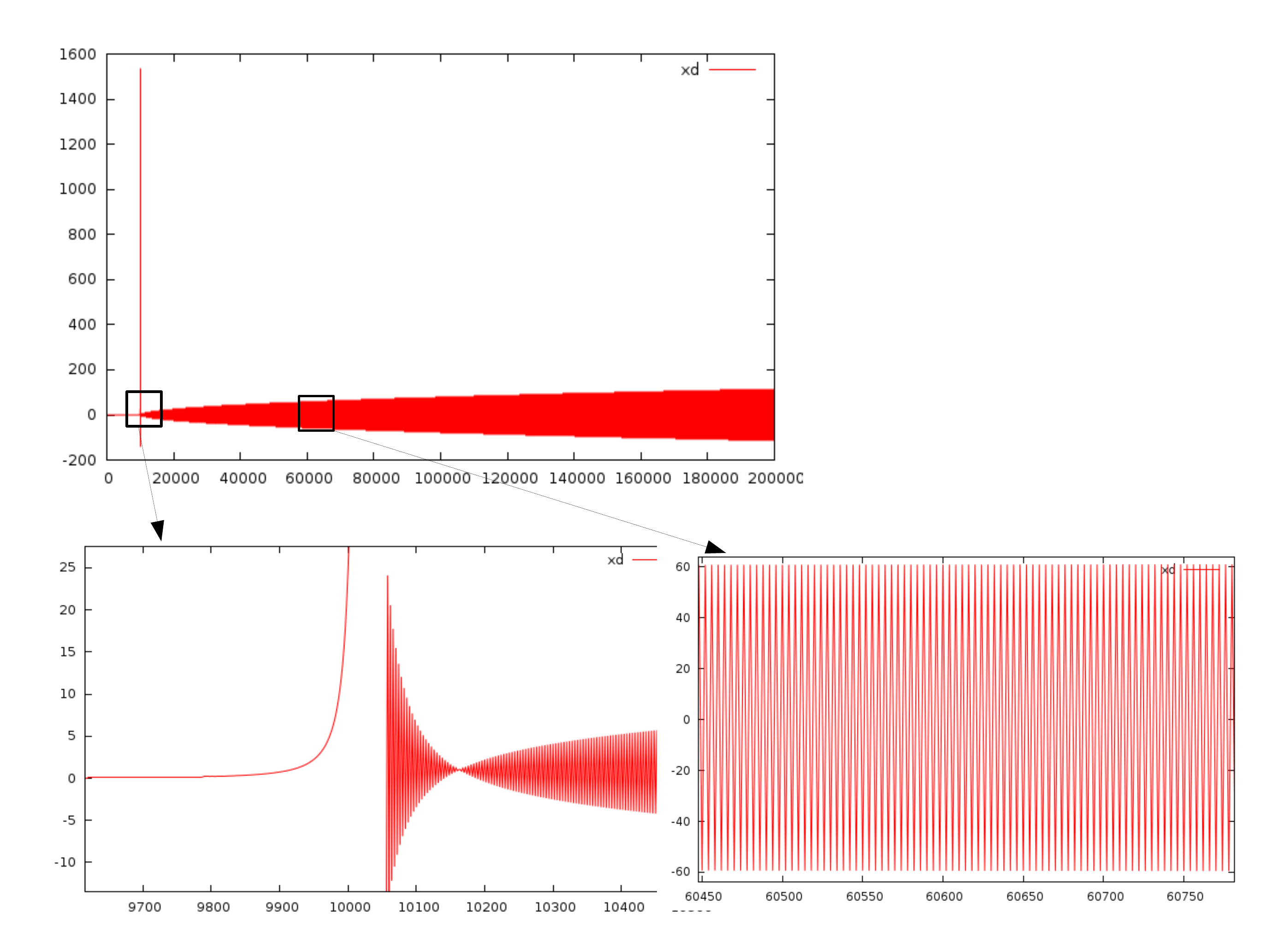}
    \caption{$dx/d\alpha$ for $\Delta t=2.0$, $\alpha=1.0$ and $\beta=1.0$ }
    \label{fig4}
\end{figure}

\begin{figure}
  \centering
    \includegraphics[scale=0.7]{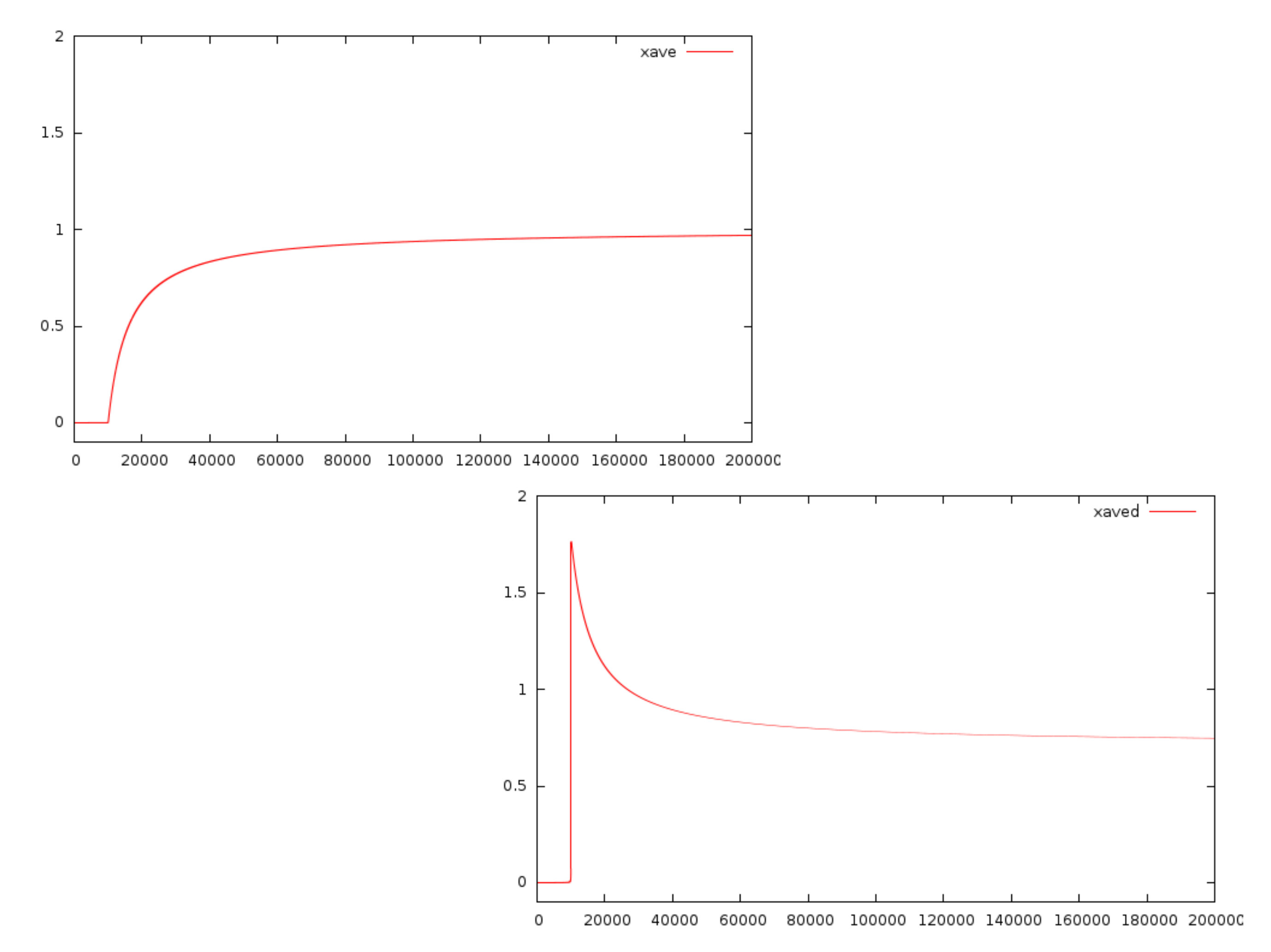}
    \caption{$J$ and $dJ/d\alpha$ for $\Delta t=2.0$, $\alpha=1.0$ and $\beta=1.0$ }
    \label{fig5}
\end{figure}

A more interesting case occurs when we take $\Delta t=2.0$. In this case (Case II), the solution is highly oscillatory as it can be seen in Fig. \ref{fig3}. It can be also observed that the oscillations tend to die out but this happens very slowly. In Fig. \ref{fig4}, the derivative of $x$ is shown. The values are also oscillatory but in contrast to the ODE solution, they tend to increase. In Fig. \ref{fig5}, the values of the object function and its derivative are shown. The averaging operation smooths the solution and removes the oscillations. Therefore, the objective function is monotone.

\begin{figure}
  \centering
    \includegraphics[scale=0.7]{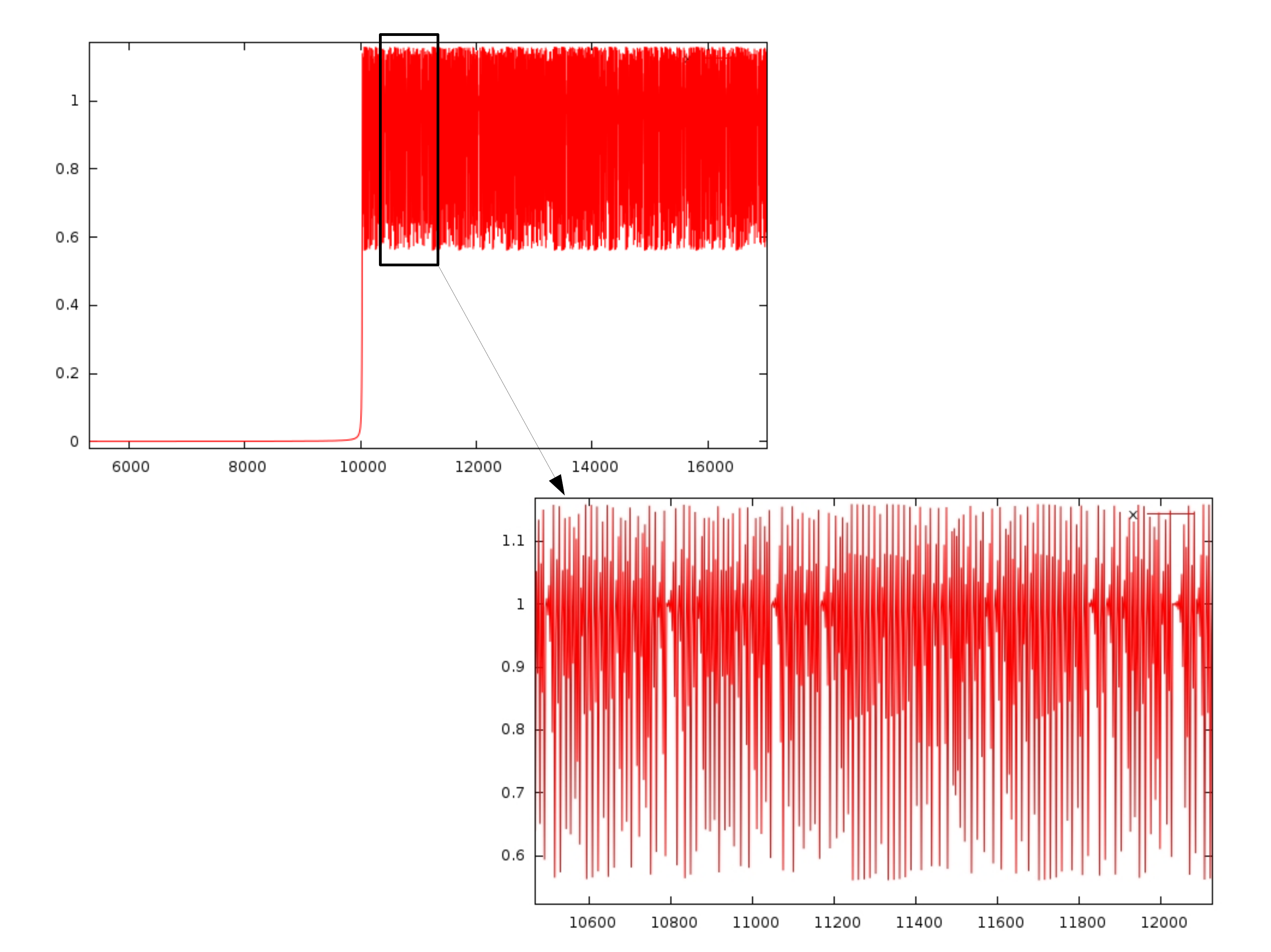}
    \caption{$x$ for $\Delta t=2.8$, $\alpha=1.0$ and $\beta=1.0$ }
    \label{fig6}
\end{figure}

\begin{figure}
  \centering
    \includegraphics[scale=0.7]{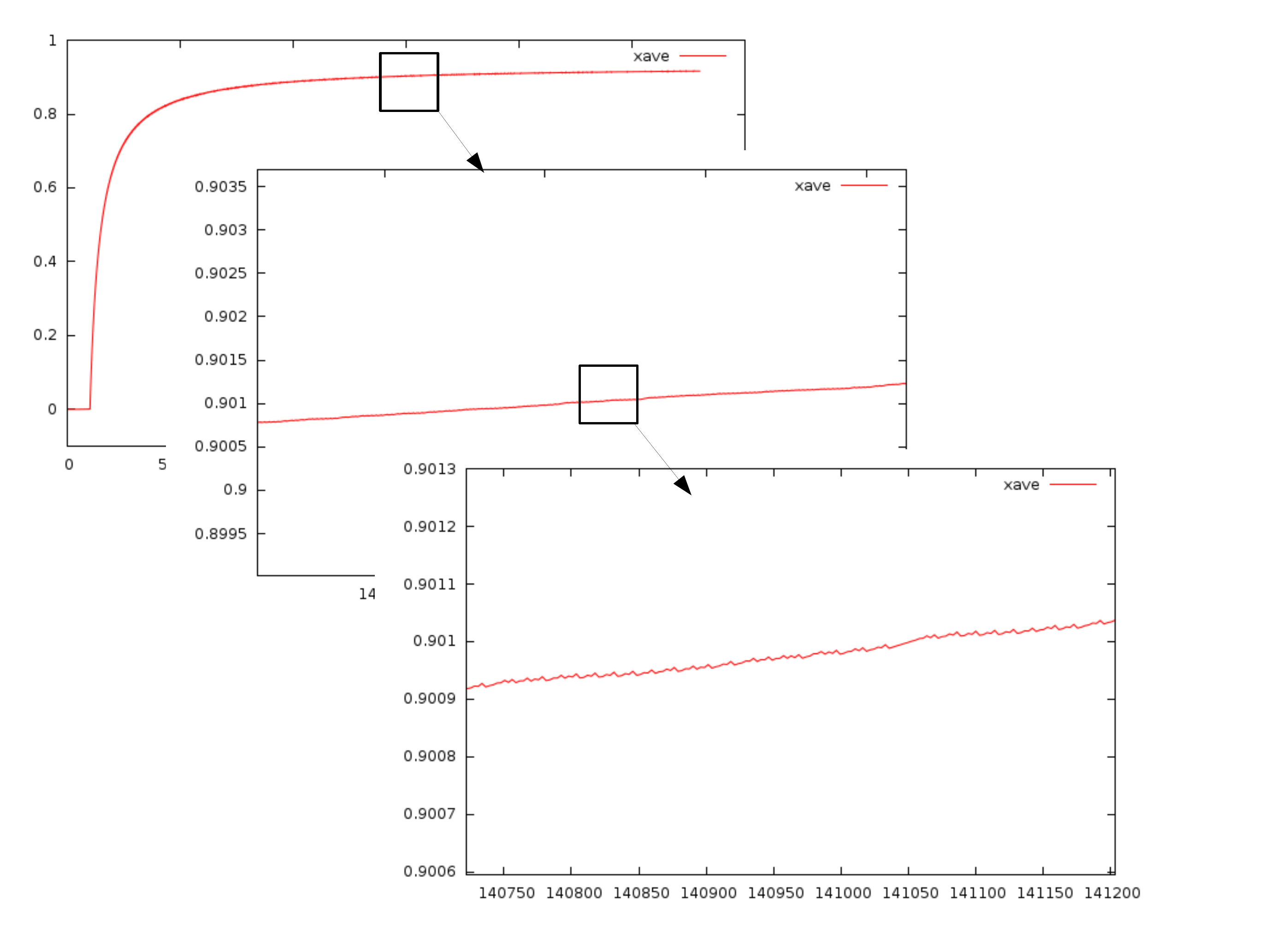}
    \caption{$J=xave$ for $\Delta t=2.8$, $\alpha=1.0$ and $\beta=1.0$ }
    \label{fig7}
\end{figure}

\begin{figure}
  \centering
    \includegraphics[scale=0.7]{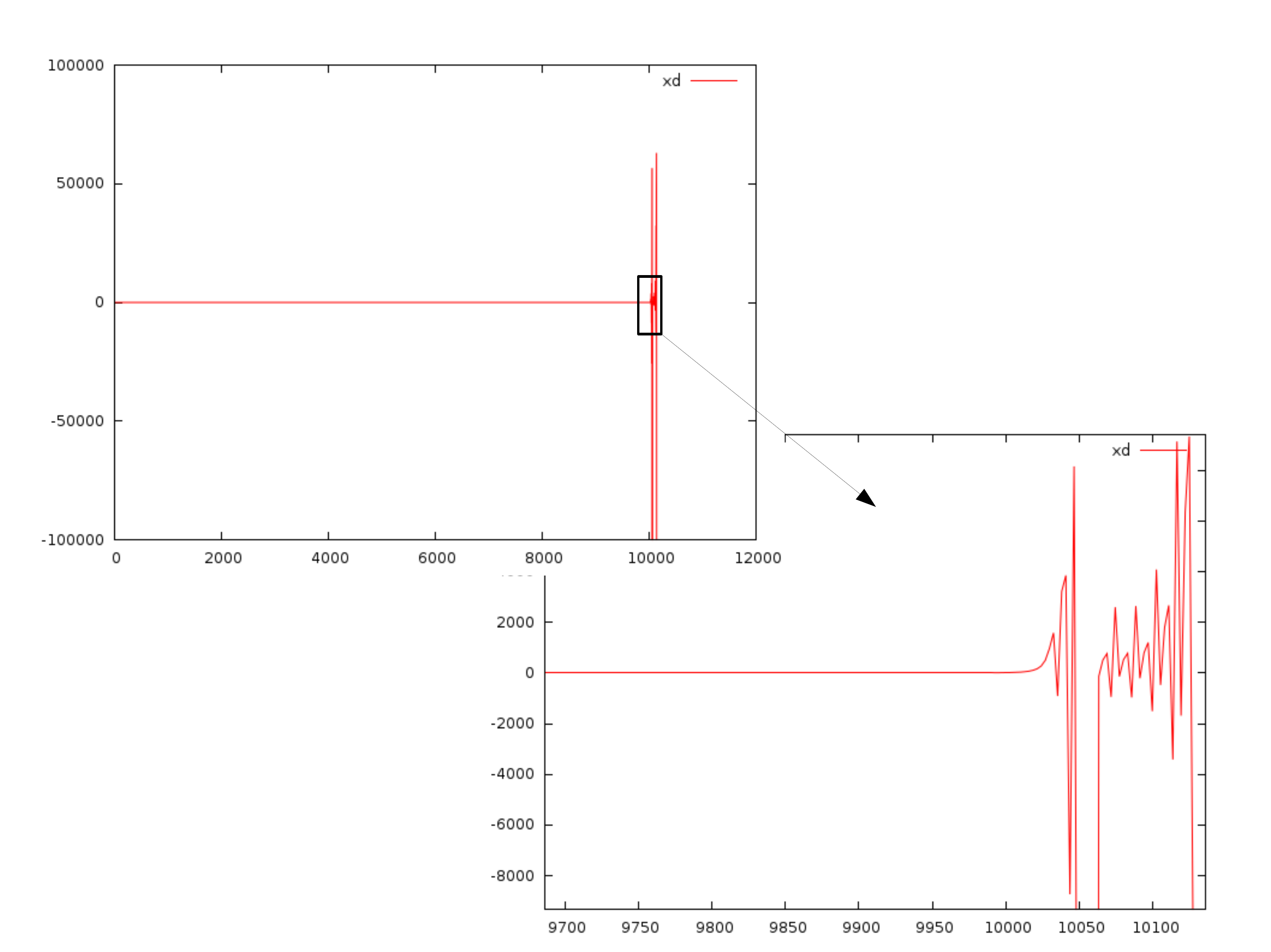}
    \caption{$xd$ for $\Delta t=2.8$, $\alpha=1.0$ and $\beta=1.0$ (values shown only up to $T=10323$ due to overflow) }
    \label{fig8}
\end{figure}

\begin{figure}
  \centering
    \includegraphics[trim=1.5cm 3.5cm 0.5cm 3cm, clip=true, scale=0.7]{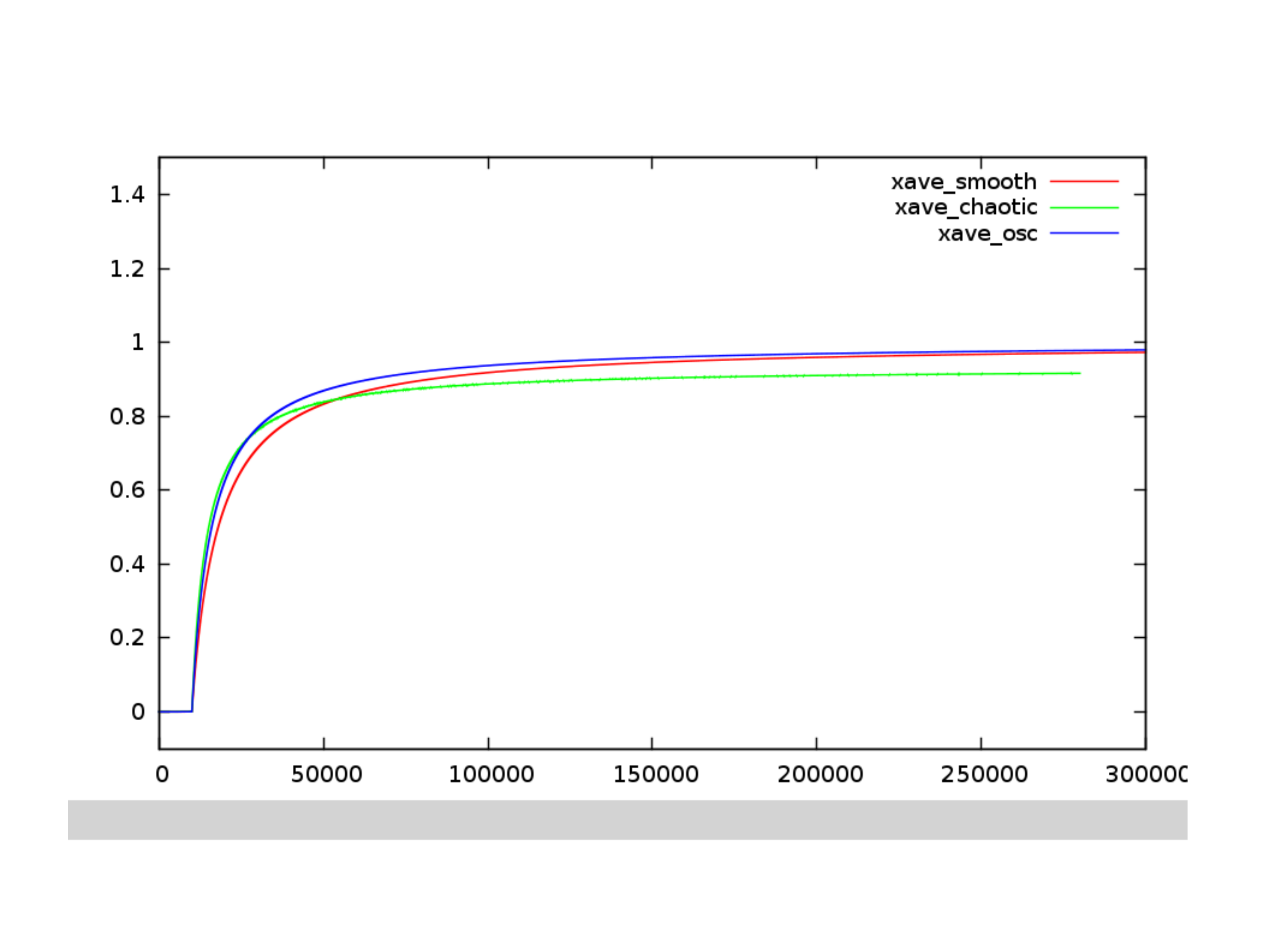}
    \caption{$J=xave$ for $\Delta t=1.0,2.0,2.8$, $\alpha=1.0$ and $\beta=1.0$ }
    \label{fig9}
\end{figure}

Finally, the most interesting case occurs when $\Delta t$ is neither large nor small. Certain values of $\Delta t$ leads to a chaotic solution. In Fig. \ref{fig6}, the solution $x$ in such a case is shown with $\Delta t=2.8$. As it can be seen from the figure, the solution becomes chaotic after a certain time step. If we look at the time averaged curve in Fig. \ref{fig7}, we see a different picture. At the first glance, the objective function looks monotone and smooth. If we zoom further into the picture, however, it can be observed that the solution has some wiggles. These wiggles are also distributed in a chaotic, random-like manner. In Fig. \ref{fig8}, the derivative values are shown. Up to a certain time-step, values grow slowly. Immediately after the primal solution becomes chaotic, derivatives grow very rapidly in an oscillating manner finally leading to an overflow. The derivatives of the objective function behave also in a similar way. In fact this result is not surprising, since the local solution at a certain time iteration becomes totally unpredictable. On the other hand, if we just look at the time-averaged quantities it is hard to notice the chaos. In Fig. \ref{fig9}, the time-averaged solution is shown again for three different values of $\Delta t$. Qualitative there is very little difference between the results and the steady-state values are pretty close to each other. 

In conclusion, the presented example shows how a purely deterministic system becomes chaotic for some numerical setup. Although ODE solutions can be within tolerable limits, differentiating chaotic solution causes the complete failure of the tangent linear solver. Since the adjoint calculus is mathematically same as the tangent linear mode operations, we can also expect the same behavior for the adjoint codes. Obviously, the remedy is to just get away from the unstable zone either by choosing a smaller $\Delta t$ or using a different integrating scheme. 

\bibliographystyle{unsrt}
\bibliography{full-paper}

\end{document}